\documentclass[a4papersize]{article}

\usepackage{amsmath,amsthm, amsxtra,amssymb,latexsym, amscd}
\usepackage[mathscr]{eucal}
\newtheorem{theorem}{Theorem}[section]
\newtheorem{corollary}[theorem]{Corollary}
\newtheorem{lemma}[theorem]{Lemma}

\newtheorem{example}[theorem]{Example}

\newtheorem{definition}[theorem]{Definition}
\newtheorem{remark}[theorem]{Remark}

\DeclareMathOperator{\depth}{depth}

\DeclareMathOperator{\Hom}{Hom}
\DeclareMathOperator{\Ann}{Ann}

\DeclareMathOperator{\Ass}{Ass}
\DeclareMathOperator{\Supp}{Supp}

\DeclareMathOperator{\Var}{Var}
\DeclareMathOperator{\nCM}{nCM}

\DeclareMathOperator{\Att}{Att}

\DeclareMathOperator{\Spec}{Spec}
\DeclareMathOperator{\p}{\frak p}
\DeclareMathOperator{\q}{\frak q}
\DeclareMathOperator{\m}{\frak m}

\begin{document}
\large
\centerline{\Large {\bf  ON SEQUENTIALLY COHEN-MACAULAY MODULES AND   }}
\smallskip
\centerline{\Large {\bf SEQUENTIALLY GENERALIZED COHEN-MACAULAY MODULES}}
\smallskip

\medskip
\vskip 0.7cm
\centerline { NGUYEN XUAN LINH}
\centerline {  Hanoi  University of Civil Engineering}
\centerline{Hanoi, Vietnam}
\centerline {e-mail: linhnx@huce.edu.vn}

\vskip 0.4cm

\centerline {LE THANH NHAN}
\centerline { Ministry of Education and Training, Hanoi, Vietnam}
\centerline {and Thai Nguyen University of Sciences, Thai Nguyen, Vietnam}
\centerline {E-mail: nhanlt2014@gmail.com}

\vskip 0.4cm
\smallskip
\vskip 1cm

\noindent{\bf Abstract} {\footnote{ {\it{Key words and phrases: }} Sequentially Cohen-Macaulay module; Sequentially generalized Cohen-Macaulay module; Sequential sequence; Sequential f-sequence; Local cohomology module. \hfill\break
  {\it{2000 Subject  Classification: }} 13E05, 13C14, 13D45.  \hfill\break
   {The first author is funded by Hanoi University of Civil Engineering (HUCE) under grant number 20-2022/KHXD-TD.}\\
   The second author is supported by the Vietnam National Foundation for Science and Technology Development (Nafosted) under grant number 101.04-2023.31.}.
     We introduce  the notions of {\it sequential  sequence} and {\it sequential f-sequence} in order to characterize sequentially Cohen-Macaulay modules and sequentially generalized Cohen-Macaulay modules. Let $R$ be a Noetherian local ring and $M$ a finitely generated $R$-module. We show that $M$ is sequentially Cohen-Macaulay (resp. sequentially generalized Cohen-Macaulay) if and only if  there exists a system of parameters of $M$ that is an $M$-sequential sequence (resp. each generalized regular sequence s.o.p of $M$ is an $M$-sequential f-sequence) and $R/\Ann_R(M)$ is a quotient of a Cohen-Macaulay local ring. As an application, we give new characterizations of Cohen-Macaulay modules and generalized Cohen-Macaulay modules. 

\section{Introduction}

\ \ \ \ Throughout this paper, let $(R,\frak m )$ be a Noetherian local ring and  $M$ a finitely generated $R$-module with $\dim_R(M)=d$.  

  It is well-known that $M$ is Cohen-Macaulay if and only if there exists a system of parameters (s.o.p for short) of $M$ that is a regular sequence, if and only if each s.o.p of $M$ is a regular sequence.   A  natural extension of the notion of regular sequence is that of filter regular sequence (see Definition \ref{D:2}) introduced by  Cuong-Schenzel-Trung \cite{CST}.   Note that  $M$ is generalized Cohen-Macaulay if and only if  $R/\Ann_R(M)$ is a quotient of a Cohen-Macaulay local ring and each s.o.p of $M$ is a filter regular sequence (see \cite[Hilfssatz 3.4, Satz 3.8]{CST} and Lemma \ref {L:2bs}).

An important generalization of the notion of Cohen-Macaulay module is that of sequentially Cohen-Macaulay module introduced almost at the same time  by  R. Stanley \cite{St}   in the graded setting  and by P. Schenzel \cite{Sch} in the local setting. The notion of sequentially generalized Cohen-Macaulay module was defined in \cite{CN} in a natural way. Let $H^0_{\frak m}(M)=D_t\subset \ldots \subset D_1\subset D_0=M$ be the dimension filtration of $M$, it means that  each $D_i$ is the largest submodule of $M$ of dimension less than $\dim_R(D_{i-1})$, see \cite{Sch}, \cite{CN}.  We say that  $M$  is sequentially Cohen-Macaulay (resp. sequentially generalized Cohen-Macaulay)  if each quotient module $D_{i-1}/D_i$ is  Cohen-Macaulay (resp. generalized Cohen-Macaulay). The sequential  Cohen-Macaulayness and sequential generalized Cohen-Macaulayness were studied from many different aspects, see \cite{St, Sch, HS, CN, TY, CC07, GHS, G, NDC}. We can see \cite{CSSS} for  a condensed summary about sequentially Cohen-Macaulay modules.

The notion of strict f-sequence was introduced in \cite{CMN} in order to study the finiteness of certain sets of prime ideals and the length of generalized fractions. This notion is a useful tool for solving some problems about canonical modules, structure of modules, preservation under small perturbation, see \cite{ANT, BN, MQ, NDC}. In this paper, we introduce  the notions of sequential  sequence and sequential f-sequence, which are slightly modified from the notion of strict f-sequence, in order to establish an analogue for the sequentially Cohen-Macaulay modules and sequentially generalized Cohen-Macaulay modules of the
above well-known characterizations of Cohen-Macaulay modules and generalized Cohen-Macaulay modules, where the roles of regular sequences and filter regular sequences are respectively replaced by that of sequential  sequences and sequential f-sequences. 

 Following I. G. Macdonald \cite{Mac}, the set of  attached primes of  an Artinian $R$-module $A$  is denoted by $\Att_RA.$ Note that the  local cohomology module $H^i_{\frak m}(M)$ is Artinian for any  $i.$

\begin{definition} \label{D:1} {\rm (a) An element $x\in \frak m$ is said to be an {\it $M$-sequential element} if $x\notin\frak p$ for all prime ideals $\frak p\in\bigcup_{j=1}^d\Att_RH^j_{\frak m}(M)$. A sequence  $x_1, \ldots , x_t$ of elements in $\frak m$ is said to be an {\it $M$-sequential sequence}  if $x_i$ is an $M/(x_1, \ldots , x_{i-1})M$-sequential element for all $i\leq t$.

 (b) An element $x\in \frak m$ is said to be an {\it $M$-sequential f-element}  if $x\notin\frak p$ for all prime ideals $\frak p\in\bigcup_{j=2}^d\Att_RH^j_{\frak m}(M)\setminus\{\frak m\}$. A sequence $x_1, \ldots , x_t$ of elements in $\frak m$ is said to be an {\it $M$-sequential f-sequence} if $x_i$ is $M/(x_1, \ldots , x_{i-1})M$-sequential f-element for all $i\leq t$}.
\end{definition}

The following theorem is the first main result of this paper, which gives characterizations of sequentially Cohen-Macaulay modules in term of sequential sequence. We note that each $M$-sequential sequence is a filter regular sequence (see Remark  \ref{R:1}(b)), and this theorem shows that if $M$ is sequentially Cohen-Macaulay then the converse statement is true.

\begin{theorem} \label{T:1a}  The following statements are equivalent:
\begin{itemize}
\item[\rm{(a)}] $M$ is sequentially Cohen-Macaulay.
\item[\rm{(b)}]  $R/ \Ann_R(M)$ is a quotient of a Cohen-Macaulay local ring and there exists a s.o.p of $M$ that is an $M$-sequential sequence.
\item[\rm{(c)}] $R/ \Ann_R(M)$ is a quotient of a Cohen-Macaulay local ring and  each filter regular sequence s.o.p of $M$ is an $M$-sequential sequence.
\end{itemize}
\end{theorem}

The notion of generalized regular sequence (see Definition \ref{D:2}) is a natural extension of the more known notions of regular sequence and filter regular sequence, which was introduced in \cite{Nh} in order to study the finiteness of associated primes of local cohomology modules.  The detail relations among regular sequences, filter regular sequences, generalized regular sequences, sequential sequences, sequential f-sequences will be clarified in Remark \ref{R:1}.

The following theorem is the second main result of this paper, which gives a characterization of sequentially generalized Cohen-Macaulay modules in term of  sequential f-sequence. We note that each $M$-sequential f-sequence is a generalized regular sequence (see Remark \ref{R:1}(b)), and this theorem shows that if $M$ is sequentially generalized Cohen-Macaulay then the converse statement is true.  
 
\begin{theorem} \label{T:2a}  The following statements are equivalent:
\begin{itemize}
\item[\rm{(a)}] $M$ is sequentially generalized Cohen-Macaulay.
\item[\rm{(b)}] $R/\Ann_R(M)$ is a quotient of a Cohen-Macaulay local ring and each generalized regular sequence s.o.p of $M$ is an $M$-sequential f-sequence. 
\end{itemize}
\end{theorem}

 In the next section, we present some preliminaries that will be used in the sequel.  In Section 3, we prove Theorem \ref{T:1a} and Theorem \ref{T:2a}.  As an application, we obtain new characterizations of Cohen-Macaulay modules and generalized Cohen-Macaulay modules, see Corollary \ref{C:3}. The key techniques for proving the main results  are the shifted principles for attached primes of local cohomology modules \cite{NQ} and a measure for non sequential Cohen-Macaulayness of modules \cite{NDC}.

\section{Preliminaries}

\ \ \ \ We always assume that $(R,\frak m )$ is a Noetherian local ring, $M$ is a finitely generated $R$-module with $\dim_R(M)=d$. Let $\widehat R$ and $\widehat M$ respectively denote the $\frak m$-adic completion of $R$ and $M$. For an ideal $I$ of $R$, denote by $\Var (I)$ the set of all prime ideals of $R$ containing $I$.

The set of attached primes introduced by I. G. Macdonald  \cite{Mac} makes an important role in the study of Artinian modules, which is  similar to the set of associated primes in the study of Noetherian modules.  Let $A$ be an Artinian $R$-module. Then $A$ has a minimal secondary representation $A=A_1+\ldots +A_n,$ where each $A_i$ is $\p_i$-secondary. The set $\{\p_1,\ldots ,\p_n\}$ is independent of the choice of the minimal  secondary representation of $A$. This set is denoted by  $\Att_RA$ and it is called {\it the set of attached primes} of $A$.     Here are some basic properties of attached primes of Artinian modules, see \cite{Mac}, \cite[8.2.4, 8.2.5]{BS}.

\begin{lemma}  \label{L:1a} Let $A$ be an Artinian $R$-module. Then
\begin{itemize}
\item[\rm{(a)}] $\min\Att_RA=\min\Var (\Ann_RA).$ In particular,  $A\neq 0$ if and only if $\Att_RA\neq \emptyset$.
\item[\rm{(b)}] $A$  has a natural structure as an Artinian $\widehat R$-module and $\Att_RA=\{\frak P\cap R\mid \frak P\in\Att_{\widehat R}A\}.$
\end{itemize}
\end{lemma}

Denote by $H^i_{\frak m}(M)$ the $i$-th local cohomology module of $M$ with respect to $\frak m$.  Note that   $H^i_{\frak m}(M)$ is  an Artinian $R$-module for all integers  $i\geq 0$. Here are some properties of attached primes of $H^i_{\frak m}(M)$, see \cite[11.3.9, 7.3.2]{BS},  \cite[Lemma 2.2]{NDC}.  

\begin{lemma}   \label{L:1b} The following statements are true.
\begin{itemize}
\item[\rm{(a)}] If $\frak p\in\Ass_R(M)$ and $\dim (R/\frak p)=i$ then $\frak p\in\Att_RH^i_{\frak m}(M).$ If  $\frak p\in\Att_RH^i_{\frak m}(M)$ with  $\dim (R/\frak p)=i$ and $R$ is a quotient of a Cohen-Macaulay local ring then $\frak p\in\Ass_R(M)$.
 \item[\rm{(b)}] If $R$ is a quotient of a Cohen-Macaulay ring and $\frak p\in\Att_RH^i_{\frak m}(M)$ then $\dim(R/\frak p)\leq i$. 
\item[\rm{(c)}] $\Att_RH^d_{\frak m}(M)=\{\frak p \in \Ass_R(M) \mid \dim (R/\frak p)=d\}.$
\end{itemize}
\end{lemma}

The following shifted principles for attached primes of local cohomology modules under localization and completion were established in \cite[Theorem 1.1]{NQ}. These shifted principles will be used as one of key techniques for the proof of the main results of this paper.

\begin{lemma}\label{L:NQ} The following statements are equivalent:
\begin{itemize}
\item[\rm{(a)}] $R$ is a quotient of a Cohen-Macaulay local ring.
\item[\rm{(b)}] $\Att_{R_{\frak p}}H^{i-\dim(R/\frak p)}_{\frak p R_{\frak p}}(M_{\frak p})=\{\frak q R_{\frak p}\mid \frak q\in\Att_RH^i_{\frak m}(M), \frak q\subseteq \frak p\}$ for all  $M$ and all $i$.
\item[\rm{(c)}] $\Att_{\widehat R}H^i_{\frak m}(M)=\underset{\frak p\in\Att_RH^i_{\frak m}(M)}{\bigcup}\Ass(\widehat R/\frak p\widehat R)$ for all $M$ and all $i$.
\end{itemize}
\end{lemma}
For convenience, we recall the notion of filter regular sequence introduced in \cite{CST} and the notion of generalized regular sequence introduced in \cite{Nh}.  

\begin{definition} \label{D:2} {\rm (a) An element $x\in\frak m$ is said to be a {\it filter regular element} (f-element for short) of $M$ if $x\notin\frak p$ for all $\p\in \Ass_R(M)$ satisfying $\dim (R/\frak p)>0$.  A sequence $x_1, \ldots , x_t$ of elements in $\frak m$ is called a {\it filter regular sequence} (f-sequence for short) of $M$ if  $x_i$ is an f-element  of $M/(x_1, \ldots , x_{i-1})M$ for all $i\leq t$.

(b) An element $x\in\frak m$ is  a {\it generalized regular element} of $M$ if $x\notin\frak p$ for all $\p\in \Ass_R(M)$ satisfying $\dim (R/\frak p)>1$. A sequence $x_1, \ldots , x_t$ of elements in $\frak m$ is said to be a {\it generalized regular sequence} if $x_i$ is a generalized regular element of $M/(x_1, \ldots , x_{i-1})M$ for all $i\leq t$.}
\end{definition}

By Prime Avoidance, for any integer $t>0,$ there always exists an f-sequence (hence a generalized regular sequence) of $M$ of length $t.$ It is clear that each f-sequence  of length $t\leq d$ is a part of s.o.p of $M$; each generalized regular sequence  of length $t\leq d-1$ is a part of s.o.p of $M$.   The following remark gives  relations among the notions of regular sequence, f-sequence, generalized regular sequence, sequential sequence, sequential f-sequence.

\begin{remark} \label{R:1} {\rm (a) Each regular sequence  is an f-sequence;  and each f-sequence  is a generalized regular sequence. The converse statements are not true. By \cite[Example 5.7]{Nh}, for given integers $d\geq 5$ and $0\leq i<j<k\leq d-3$, there exist a finitely generated $R$-module $M$ and an ideal $I$ of $R$ such that $i, j,  k$ are repectively the common lengths of maximal regular sequences, maximal f-sequences, maximal generalized regular sequences of $M$ in $I$.   

(b)  It follows by Lemma \ref{L:1b}(a)  that  $\Ass_R(M)\subseteq \bigcup_{j=0}^{\dim_R (M)}\Att_RH^j_{\frak m}(M)$ for all finitely generated $R$-module $M$. Therefore, each sequential sequence is an f-sequence;   each sequential f-sequence  is a generalized regular sequence.

(c) Regular sequences are not necessarily sequential sequences. Moreover, regular sequences are not necessarily sequential f-sequences.  For example, let $R=K[[x,y,z,t,w]]$ be the formal power series ring over a field $K$ and $M=R/(x,y)\cap \frak (z,t)$. Set $\frak m=(x,y,z,t,w)$, $\frak p_1=(x,y)$, $\frak p_2=(z,t)$, $\frak p_3=(x,y,z,t)$. Then $\Att_RH^0_{\frak m}(M)=\Att_RH^1_{\frak m}(M)=\emptyset$, $\Att_R H^2_{\frak m}(M)=\{\frak p_3\}$ and $\Att_R H^3_{\frak m}(M)=\{\frak p_1, \frak p_2\}$. Therefore, $x+z$ is an $M$-regular element but it is not a sequential f-element.

(d) Sequential sequences are not necessarily regular sequences; sequential f-sequences are not necessarily f-sequences. For example, let $R=K[[x, y, z]]$  be the formal power series ring over a field $K$. Let $M=R/I$, where $I=(x)\cap (y,z)\cap (x^2, y^2, z)$. Set $\frak m=(x, y, z)$, $\frak p_1=(x)$, $\frak p_2=(y,z).$ We have  $\Att_RH^2_{\frak m}(M)=\{\frak p_1\}$, $\Att_RH^1_{\frak m}(M)=\{\frak p_2\}$. Therefore, $x+z$ is an $M$-sequential element but it is not a regular element; $y$ is an $M$-sequential f-element but it is not an f-element.

(e) By weak general shifted localization principle \cite[11.3.8]{BS}, if $x_1, \ldots , x_t$ is an $M$-sequential f-sequence then $x^*_1, \ldots , x^*_t$ is an $M_{\frak p}$-sequential sequence for all $\frak p\in\Spec(R)\setminus\{\frak m\}$ containing $x_1, \ldots , x_t$. Here $x^*_i$ is the image of $x_i$ in $R_{\frak p}$ for all $i$.}
\end{remark}

The notion of dimension filtration  was first defined by P. Schenzel  \cite{Sch}. Then the definition was changed in \cite{CN}  by removing the repeated components.    

\begin{definition} {\rm A filtration  $H^0_{\frak m}(M)=D_t\subset \ldots \subset D_1\subset D_0=M$  of submodules of $M$ is said to be the {\it dimension filtration} of $M$ if   $D_i$ is the largest submodule of $M$ of dimension less than $\dim_R(D_{i-1})$ for all $i\leq t$.}  
\end{definition}  

 Note that the dimension filtration of $M$ always exists and it is unique.  Here are some properties of the dimension filtration. Following M. Nagata \cite{Na}, $M$ is said to be {\it unmixed} if $\dim (\widehat R/\frak P)=\dim_R(M)$ for all $\frak P\in\Ass_{\widehat R}(\widehat M).$ 

\begin{lemma}   \label{L:1c} Let  $H^0_{\frak m}(M)=D_t\subset \ldots \subset D_1\subset D_0=M$ be the dimension filtration of $M.$  Set $d_i:= \dim _R(D_i)$ for all $i \leq t.$ Then
\begin{itemize}
\item[\rm{(a)}] $\Ass_R(D_i)=\{\frak p \in \Ass_R(M) \mid \dim_R(R/\frak p) \leq d_i\}$. 
\item[\rm{(b)}] $\Ass_R(M/D_i)=\{\frak p \in \Ass_R(M) \mid \dim_R(R/\frak p) > d_i\}$. 
\item[\rm{(c)}] $\Ass_R(D_i/D_{i+1})=\{\frak p \in \Ass_R(M) \mid \dim_R(R/\frak p) = d_i\}$.
In particular, we have $$\Ass_R(M)\setminus \{\frak m\}=\bigcup_{i=1}^{t} \Ass_R(D_{i-1}/D_i).$$
\item[\rm{(d)}] If $R/\frak p$ is unmixed for all $\frak p\in\Ass_R(M)$ then $H^0_{\frak m\widehat R}(\widehat M)=\widehat D_t \subset \ldots \subset \widehat D_1 \subset \widehat D_0=\widehat M$ is the dimension filtration of $\widehat M.$   
\end{itemize}
\end{lemma}
\begin{proof} 
The statements {\rm (a), (b), (c)} follow by \cite[Corollary 2.3]{Sch}. 

We prove statement {\rm (d)} by induction on $t$. The case $t=0$ is clear.  Let $t\geq 1$. Then $d>0$. Let $U_1$ be the largest submodule of $\widehat M$ of dimension less than $d$. Then  $\widehat D_1 \subseteq U_1.$ Suppose that $\widehat D_1 \neq U_1$. Let $\frak P \in \Ass_{\widehat R}({U_1}/{\widehat D_1}).$  Then $\dim ({\widehat R}/{\frak P}) \leq \dim_{\widehat R} (U_1)<d$ and $\frak P \in \Ass_{\widehat R}({\widehat M}/{\widehat D_1}).$ Set $\frak p=\frak P\cap R.$ Then $\frak p \in \Ass_{R}({M}/{D_1})$ by  \cite[Lemma 3.4]{Nh} and $\frak P\in\Ass (\widehat R/\frak p\widehat R)$ by \cite[Theorem 23.2]{Mat}. Note that   $\dim (R/\frak p)=d>0$ and $\frak p\in\Ass_R(M)$ by Lemma \ref{L:1c}(c). Since $R/\frak p$ is unmixed, $\dim(\widehat R/\frak P)=d.$  This gives a contradiction. Hence $U_1=\widehat D_1.$ Since  $H^0_{\frak m}(M)=D_t \subset \ldots \subset D_1$ is the dimension filtration of $D_1$ and $R/\frak p$ is unmixed for all $\frak p\in\Ass_R(D_1)$,  we get by induction that $H^0_{\frak m\widehat R}(\widehat M)=\widehat D_t \subset \ldots \subset \widehat D_1$ is the dimension filtration of $\widehat D_1$. 
\end{proof}

The notion of sequentially Cohen-Macaulay module was introduced by R. Stanley \cite[p. 87]{St} for graded  modules and by P. Schenzel \cite{Sch} for finitely generated modules over local rings. The notion of sequentially generalized Cohen-Macaulay module was defined in \cite{CN}. 

\begin{definition} {\rm Let $H^0_{\frak m}(M)=D_t\subset \ldots \subset D_1\subset D_0=M$ be the dimension filtration of $M$. We say that   $M$ is  {\it sequentially Cohen-Macaulay} (resp.  {\it sequentially generalized Cohen-Macaulay}) if  $D_{i-1}/D_i$ is Cohen-Macaulay (resp.  generalized Cohen-Macaulay)  for all $i\leq t$}.  
\end{definition}  

   We recall the  homological characterizations of sequentially Cohen-Macaulay modules and sequentially generalized Cohen-Macaulay modules.
 Suppose that $R$ is  a quotient of a Gorenstein local  ring $(R', \m')$ with $\dim(R')=n.$  For each integer $0\leq i\leq d$, let $K^i(M):={\rm Ext}^{n-i}_{R'}(M,R')$ denote the $i$-th deficiency module of $M$. Note that $K^i(M)$   is a finitely generated $R$-module and local duality gives the isomorphism $H^i_{\m}(M)\cong \Hom_R(K^i(M), E(R/\m)),$  where  $E(R/\m)$ is the injective hull of $R/\m$, see \cite[11.2.6]{BS}. So, $\Att_RH^i_{\m}(M)=\Ass_RK^i(M)$ by \cite[10.2.20]{BS}. Moreover,   $\dim_R(K^i(M))\leq i$ and   for any $\frak p\in\Spec(R)$ with $\dim(R/\frak p)=i$ we have $\frak p\in \Ass_R(K^i(M))$ if and only if $\frak p\in\Ass_R(M)$. A good reference for  properties of $K^i(M)$ as well as the local duality is provided in \cite{S1}.

\begin{lemma} \label{L:1d}  {\rm (See \cite[Theorem 1.2]{Sch},  \cite[Theorem 5.3]{CN}).} Suppose that $R$ is  a quotient of a Gorenstein local  ring. The following statements are true.
\begin{itemize}
\item[\rm{(a)}] $M$ is  sequentially Cohen-Macaulay  if and only if $K^i(M)=0$ or  $K^i(M)$ is a Cohen-Macaulay module of dimension $i$ for all integers $i\leq d.$
\item[\rm{(b)}] $M$ is  sequentially generalized Cohen-Macaulay  if and only if $\ell_R(K^i(M))<\infty$ or  $K^i(M)$ is a generalized Cohen-Macaulay module of dimension $i$ for all integers $i\leq d.$
\end{itemize}
\end{lemma}

    Next, we recall a measure of non Cohen-Macaulayness of modules introduced by N. T. Cuong \cite{C}. Let $\underline x= (x_1,\ldots ,x_d)$  be a s.o.p of $M$ and $\underline n =(n_1,\ldots ,n_d)$ a tuple of positive integers. Set
$$I_{M,\underline x}(\underline n ):=\ell \big(M/(x_1^{n_1},\ldots ,x_d^{n_d})M\big)-n_1\ldots n_d e(\underline x ;M),$$
 where  $e(\underline x ;M)$ is the multiplicity of $M$ with respect to $\underline x.$ In general, $I_{M,\underline x}(\underline n )$ is not a polynomial for $n_1,\ldots ,n_d \gg 0$, but it takes non-negative values and bounded above by polynomials. The least degree of all polynomials bounding above the function $I_{M,\underline x}(\underline n )$ does not depend on the choice of $\underline x$ (see \cite [Theorem 2.3]{C}). This invariant is denoted by ${\rm p}(M)$ and called the {\it polynomial type} of $M$. If we stipulate the degree of zero polynomial to be $-1$, then $M$ is Cohen-Macaulay if and only if ${\rm p}(M)=-1$, and  $M$ is generalized Cohen-Macaulay if and only if ${\rm p}(M)\leqslant 0.$ Denote by $\nCM(M):=\{\frak p\in\Spec(R)\mid M_{\frak p}\ \text{is not Cohen-Macaulay}\}$ the non-Cohen-Macaulay locus of $M$. Let $U_M(0)$ be the largest submodule of $M$ of dimension less than $d$. In general, we can compute ${\rm p}(M)$ in terms of dimension of local cohomology modules $H^i_{\frak m}(M)$, dimension of the non-Cohen-Macaulay locus $\nCM(M)$ and dimension of $U_M(0)$, see \cite{C}, \cite[Lemma 2.3]{NDC}. Note that $\nCM(M)$ is not necessarily closed under Zariski topology, but it is  closed under specialization, i.e. if $\q\in\nCM (M)$ then $\p\in\nCM (M)$  for all $\q,\p\in\Spec(R)$ with $\q\subseteq \p$. So, we can define $\dim_R(\nCM (M))$ in the usual way.  

\begin{lemma} \label{L:1e} The following statements are true.
\begin{itemize}
\item[\rm{(a)}] ${\rm p}(M)={\rm p}(\widehat M)=\underset{i\leqslant d-1}{\max}\ \dim \big(\widehat R/\Ann_{\widehat R}H^i_{\frak m}(M)\big).$ In particular, ${\rm p}(M)\leq d-1.$
\item[\rm{(b)}] ${\rm p}(M)\geq \max\ \{\dim_R(\nCM (M)), \dim_R(U_M(0))\}$. If $R$ is a quotient of a Cohen-Macaulay local ring then $\nCM (M)$ is closed under Zariski topology and $${\rm p}(M)= \max\ \{\dim_R(\nCM (M)), \dim_R(U_M(0))\}.$$
\end{itemize}
\end{lemma} 

The  {\it sequential polynomial type} ${\rm sp}(M)$ of $M$  was introduced  in \cite{NDC} in order to study the non-sequential Cohen-Macaulayness of $M$. It is defined as follows
$${\rm sp}(M):=\max\ \{{\rm p} (D_{i-1}/D_i)\mid i=1, \ldots , t\},$$
where $H^0_{\frak m}(M)=D_t \subset \ldots \subset D_1 \subset D_0=M$ is the dimension filtration of $M.$ It is clear that $M$ is sequentially Cohen-Macaulay if and only if ${\rm sp}(M)=-1$, and  $M$ is sequentially generalized Cohen-Macaulay if and only if ${\rm sp}(M)\leqslant 0.$  In general, if $R$ is a quotient of a Cohen-Macaulay local ring, then the non sequentially Cohen-Macaulay locus of $M$ is closed under Zariski topology and  ${\rm sp}(M)$ is exactly the dimension of the non sequentially Cohen-Macaulay locus of $M$, see \cite[Proposition 3.2]{NDC}.  The following homological characterization of sequential polynomial type (see \cite[Theorem 4.7]{NDC}) will be used as a key technique for the proof of the main results.

\begin{lemma} \label{L:1g}  Set $D(M)=\{\dim(R/\frak p)\mid \frak p\in\Ass_R(M)\}.$ Suppose that $R$ is a quotient of a Gorenstein local ring. Set $q_1=\underset{i\notin D(M)}{\max} \dim_R(K^i(M))$ and $q_2=\underset{i\in D(M)}{\max} {\rm p}(K^i(M)).$ Then 
$${\rm sp}(M)=\max\{q_1, q_2\}.$$
\end{lemma} 

\section{Main results}

Note that if $M$ is  sequentially  Cohen-Macaulay (resp. sequentially generalized Cohen-Macaulay)  then so is $\widehat M$, cf. \cite[Theorem 4.9]{Sch}, \cite[Theorem 3.5]{NDC}. The converse statement is not true in general, see  \cite[Example 6.1]{Sch}, \cite[Example 3.7]{NDC}. 

\begin{lemma}  \label{L:2b} The following satements are true.
\begin{itemize}
\item[\rm{(a)}] $M$ is  sequentially  Cohen-Macaulay  if and only if $\widehat M$ is  sequentially  Cohen-Macaulay  and $R/\frak p$ is unmixed  for all $\frak p \in \Ass_R(M)$.
\item[\rm{(b)}] $M$ is  sequentially  generalized Cohen-Macaulay  if and only if $\widehat M$ is sequentially  generalized Cohen-Macaulay and $R/\frak p$ is unmixed  for all $\frak p \in \Ass_R(M)$.
\end{itemize}
\end{lemma}

\begin{proof}  The statement {\rm (a)} follows  by \cite[Proposition 4.6, Theorem 4.9]{Sch} and Lemma \ref{L:1c}(d). 

{\rm (b)} Suppose that $M$ is  sequentially  generalized Cohen-Macaulay. Let $\frak p\in\Ass_R(M)$ and $\frak P\in\Ass (\widehat R/\frak p\widehat R).$ If $\frak p=\frak m$ then clearly $R/\frak p$ is unmixed. So, we assume that $\frak p\neq \frak m$.     Let  $H^0_{\frak m}(M)=D_t \subset \ldots \subset D_1 \subset D_0=M$ be the dimension filtration of $M.$ Then by Lemma \ref{L:1c}(c),  $\frak p\in\Ass_R(D_{i-1}/D_i)$ and $\dim(R/\frak p)=\dim_R(D_{i-1})$ for some  $1\leq i\leq t$.  So,  $\frak P\in \Ass_{\widehat  R}(\widehat D_{i-1}/\widehat D_i)$ by \cite[Theorem 23.2]{Mat}.  Since $D_{i-1}/D_i$ is generalized Cohen-Macaulay,  so is $\widehat D_{i-1}/\widehat D_i$.  Because $\frak P\cap R=\frak p\neq \frak m$, we have  $\frak P\neq\frak m\widehat R$. It follows that  $\dim (\widehat R/\frak P)=\dim_{\widehat  R} (\widehat D_{i-1}/\widehat D_i)=\dim_R(D_{i-1}).$ Hence $R/\frak p$ is unmixed. Now the result follows by \cite[Theorem 3.5]{NDC} and Lemma \ref{L:1c}(d).
\end{proof}

It turns out by Lemma \ref{L:2b} that the sequential Cohen-Macaulayness as well as sequential generalized Cohen-Macaulayness are better described over a complete local ring. Then local duality is available, so  we can use deficiency modules instead of local cohomology modules.   This idea will be used in the proofs of the main results in the sequel.  In the following lemma we  see that  the notions of sequential element and sequential f-element are closedly related to the concept of  strongly filter regular element  defined in \cite[Definition 2.5]{S2}.

\begin{lemma} \label{L:2bbs}  (a) An element $x\in\frak m$ is an $M$-sequential elelement (resp. $M$-sequential f-element)   if and only if $x$ is $K^i(\widehat M)$-regular for  all $i\geq 1$ with  $K^i(\widehat M)\neq 0$ (resp. $x$ is a f-element of $K^i(\widehat M)$ for  all $i\geq 2$). 

(b) A sequence  $x_1, \ldots , x_t\in\frak m$   is an $M$-sequential sequence (resp.  $M$-sequential f-sequence)  if and only if it is an $\widehat M$-sequential sequence (resp.  $\widehat M$-sequential f-sequence).
\end{lemma}
\begin{proof}  It follows by local duality  and by \cite[10.2.20]{BS} that
$$\Att_{\widehat R} H^i_{\frak m}(M)=\Att_{\widehat R} H^i_{\frak m\widehat R}(\widehat M)=\Ass_{\widehat R}(K^i({\widehat M}))$$
for all $i\geq 0.$ Now the result follows by Lemma \ref{L:1a}(b). 
\end{proof}

We need recall the notion of p-standard s.o.p introduced by N. T. Cuong \cite{C}. This notion makes an important role in the study of the singularity of Cohen-Macaulay type of Noetherian rings and modules. T. Kawasaki \cite{Kaw} used p-standard s.o.p to show that a Noetherian local ring has an arithmetic Macaulayfication if and only if it is unmixed and all its formal fibers are Cohen-Macaulay. Kawasaki's theorem was extended for modules in \cite{CC17}, then it was constructed for idealization in \cite{CNN24}.    Set $\frak a(M):=\frak a_0(M)\frak a_1(M)\ldots \frak a_{d-1}(M),$ where $\frak a_i(M)=\Ann_RH^i_{\frak m}(M)$ for  $i\leq d$. An s.o.p $(x_1, \ldots , x_d)$ of $M$ is said to be a {\it p-standard s.o.p} if $x_d\in\frak a(M)$ and $x_i\in\frak a(M/(x_{i+1}, \ldots , x_d)M)$ for all $i<d.$ 
 Note that if  $(x_1,\ldots, x_d)$ a p-standard s.o.p of $M$ then there exist   $\lambda _0,\ldots, \lambda _d\in\Bbb N$ such that 
$$\ell(M/(x_1^{n_1}, \ldots, x_d^{n_d})M)=\sum\limits_{i=0}^d\lambda _i n_1\ldots n_i$$ for all integers $n_1,\ldots, n_d \geq 1,$ see \cite[Theorem 2.6]{C}.  In general, p-standard s.o.p does not necessarily exist. By \cite[Theorem 1.2, 1.3]{CC17},  $M$ admits a p-standard s.o.p if and only if $R/\Ann_R(M)$ is a quotient of a Cohen-Macaulay local ring.

\begin{lemma} \label{L:2bs} If $M$ is sequentially generalized Cohen-Macaulay then $R/\Ann_R(M)$ is a quotient of a Cohen-Macaulay local ring.
\end{lemma}

\begin{proof}  By  \cite[Theorems 1.2, 1.3]{CC17}, it is enough to show that $M$ admits a p-standard s.o.p. We prove this by induction on $d$. The case $d\leq 1$ is clear. Let $d>1.$  We first claim that if $\dim (R/\frak a_i(M))>0$ then $\dim (R/\frak a_i(M))=i$ for all integers $i\geq 1,$ where $\frak a_i(M)=\Ann_RH^i_{\frak m}(M)$. In fact,  asssume that $\dim (R/\frak a_i(M))>0.$ Let $\frak p\in\min\Var(\frak a_i(M))\setminus\{\frak m\}.$ By Lemma \ref{L:1a}(a),  $\frak p\in\Att_R(H^i_{\frak m}(M))\setminus\{\frak m\}.$ Hence, there exists $\frak P\in \Att_{\widehat R}(H^i_{\frak m}(M))\setminus \{\frak m\widehat R\}$ such that $\frak P\cap R=\frak p,$ see Lemma \ref{L:1a}(b).  So, $\frak P\in \Ass_{\widehat R}(K^i({\widehat M}))\setminus \{\frak m\widehat R\}$ by \cite[10.2.20, 11.2.6]{BS}.
 Since $M$ is sequentially generalized Cohen-Macaulay, so is $\widehat M$ by Lemma \ref{L:2b}(b). Therefore, we get by Lemma \ref{L:1d}(b) that  $K^i({\widehat M})$ is generalized Cohen-Macaulay of dimension $i$. Hence $\dim(\widehat R/\frak P)=i$ and hence $\frak P\in\Ass_{\widehat R}(\widehat M)$ by Lemma \ref{L:1b}(a). So, $\frak p\in\Ass_R(M)$ by \cite[Lemma 3.4]{Nh} and hence $\frak P\in\Ass (\widehat R/\frak p\widehat R)$ by \cite[Theorem 23.2]{Mat}. It follows by Lemma \ref{L:2b}(b) that $R/\frak p$ is unmixed. So, $\dim (R/\frak p)=\dim(\widehat R/\frak P)=i.$ Therefore, $\dim (R/\frak a_i(M))=i$ and the claim follows.  

We have by the claim  that $\dim (R/\frak a(M))<d,$ where $\frak a(M)=\frak a_0(M)\ldots \frak a_{d-1}(M).$  Let 
$$H^0_{\frak m}(M)=D_t\subset D_{t-1}\subset \ldots \subset D_1\subset D_0=M$$ be the dimension filtration of $M$.  Since $\dim_R(D_1)<d$ and $\dim (R/\frak a(M))<d,$ we can choose a parameter element $x_d$ of $M$ such that  $x_d\in\frak a(M)\cap \Ann_R(D_1).$  Since $\dim_R(0:_Mx^2_d)<d$, we have $D_1\subseteq (0:_Mx_d)\subseteq  (0:_Mx^2_d)\subseteq D_1$. So, $D_1=(0:_Mx_d)=(0:_Mx^2_d).$ Let $x_dm\in D_1\cap x_dM.$ Then $x_dm\in (0:_Mx_d)$. Hence $x_d^2m=0.$ So, $m\in  (0:_Mx^2_d).$ Hence $m\in (0:_Mx_d)$ and hence $x_dm=0.$ Therefore, $D_1\cap x_dM=0.$  So, $D_i\cap x_dM=0$ for all $i\geq 1.$  Set $\overline D_i=(D_i+x_dM)/x_dM$ for $i\geq 0$.  Then we have  $\overline D_i\cong D_i$ for all $i\geq 1.$ So, $\overline D_i/\overline D_{i+1}$ is generalized Cohen-Macaulay  and $\dim_R(\overline D_{i+1})<\dim_R(\overline D_i)$ for all $i\geq 1.$ Since $M/D_1$ is generalized Cohen-Macaulay of dimension $d$ and $x_d$ is a parameter element of $M/D_1$,  it follows that $\overline D_0/\overline D_1$ is generalized Cohen-Macaulay of dimension $d-1$.  If $\dim_R(D_1)<d-1,$ then from the filtration $\overline D_t\subset \overline D_{t-1}\subset \ldots \subset \overline D_1\subset \overline D_0=M/x_dM$, it follows by \cite[Lemma 44(iii)]{CN} that $M/x_dM$ is sequentially generalized Cohen-Macaulay. Suppose that  $\dim_R(D_1)=d-1.$ From the exact sequence $0\to \overline D_1/\overline D_2\to \overline D_0/\overline D_2\to \overline D_0/\overline D_1\to 0$ with notice that $\overline D_1/\overline D_2$ and $\overline D_0/\overline D_1$ are generalized Cohen-Macaulay of dimension $d-1$, it follows that $\ell_R(H^j_{\frak m}(\overline D_0/\overline D_2))<\infty$ for all $j<d-1.$ Hence $\overline D_0/\overline D_2$ is  generalized Cohen-Macaulay of dimension $d-1$. Therefore, from from the filtration $\overline D_t\subset \overline D_{t-1}\subset \ldots \subset \overline D_2\subset \overline D_0=M/x_dM$, it follows by \cite[Lemma 44(iii)]{CN} that $M/x_dM$ is sequentially generalized Cohen-Macaulay. By induction, there exists a p-standard s.o.p $x_1, \ldots , x_{d-1}$ of $M/x_dM.$ Hence $x_1, \ldots , x_d$ is a p-standard s.o.p of $M$.
 \end{proof}

Now we prove Theorem \ref{T:1a}.

\begin{theorem} \label{T:1}  The following statements are equivalent:
\begin{itemize}
\item[\rm{(a)}] $M$ is sequentially Cohen-Macaulay.
\item[\rm{(b)}]  $R/ \Ann_R(M)$ is a quotient of a Cohen-Macaulay local ring and there exists a s.o.p of $M$ that is an $M$-sequential sequence.
\item[\rm{(c)}] $R/ \Ann_R(M)$ is a quotient of a Cohen-Macaulay local ring and  each f-sequence s.o.p of $M$ is an $M$-sequential sequence.
\item[\rm{(d)}]  $R/ \frak p$ is unmixed for all $\frak p \in \Supp_R(M)$ and there exists a s.o.p of $M$ that is an $M$-sequential sequence.
\item[\rm{(e)}] $R/ \frak p$ is unmixed for all $\frak p \in \Supp_R(M)$ and each f-sequence s.o.p of $M$ is an $M$-sequential sequence.
\end{itemize}
\end{theorem}

\begin{proof} (a)$\Rightarrow$(e). From assumption (a), we get by \cite[Proposition 4.6(b)]{Sch} that $R/\frak p$ is unmixed for all $\frak p\in\Supp_R(M)$.  Let $(x_1, \ldots , x_d)$ be an  f-sequence s.o.p of $M$. We prove by induction on $d$ that  $(x_1, \ldots , x_d)$ is an $M$-sequential sequence. The case $d\leq 1$ is clear by Lemma \ref{L:1b}(c). Let $d\geq 2.$    Suppose in contrary that  $x_1$ is not an $M$-sequential element.  By Lemma \ref{L:2bbs}(a) that $x_1\in \frak P$ for some $\frak P\in \Ass_{\widehat R}(K^i({\widehat M}))$ with $i\geq 1.$  Hence $K^i({\widehat M})\neq 0.$ From assumption (a), we get by Lemma \ref{L:2b}(a) that  $\widehat M$ is sequentially Cohen-Macaulay.  So, we get by Lemma \ref{L:1d}(a) that $K^i({\widehat M})$ is Cohen-Macaulay of dimension $i$. Therefore $\dim (\widehat R/\frak P)=i.$  So,  $\frak P\in\Ass_{\widehat R}(\widehat M)$. Therefore $\frak p=\frak P\cap R\in\Ass_R(M)$, see \cite[Lemma 3.4]{Nh}.  Since $x_1\in\frak p$ and $x_1$ is an f-element,  $\frak p=\frak m$. Hence $\frak P=\frak m\widehat R$, and hence  $1\leq i=\dim (\widehat R/\frak P)=0.$ This gives a contradiction, so $x_1$ is  an $M$-sequential element.  From assumption (a) with notice that $x_1$ is an f-element of $M$,  we get  by  \cite[Lemma 4.1(i)]{NDC}  that $M/x_1M$ is sequentially Cohen-Macaulay. Since $(x_2, \ldots , x_d)$ is an f-sequence s.o.p of $M/x_1M$,  we have by induction that   $(x_2, \ldots , x_d)$  is an $M/x_1M$-sequential sequence. Hence $(x_1, \ldots , x_d)$  is an $M$-sequential sequence.

The statement (e)$\Rightarrow$(d) is obvious.

(d)$\Rightarrow$(a). By assumption (d), we can take a s.o.p $(x_1, \ldots , x_d)$ of $M$ that is an $M$-sequential sequence. By Lemma \ref{L:2b}(a) and the unmixedness of $R/\frak p$ for all $\frak p \in \Ass_R(M)$, it is enough to prove that $\widehat M$ is sequentially Cohen-Macaulay. We prove this by induction on $d$. The case $d\leq 1$ is obvious. Let $d\geq 2.$  Note that  $\dim_R(M/x_1M)=d-1$, $R/\frak p$ is unmixed for all $\frak p\in\Supp_R(M/x_1M)$ and $(x_2, \ldots , x_d)$ is a s.o.p of $M/x_1M$ which is an $M/x_1M$-sequential sequence. 
So we have by induction that $\widehat M/x_1\widehat M$ is sequentially Cohen-Macaulay. By Lemma \ref{L:1d}(a), $K^i({\widehat M/x_1\widehat M})=0$ or $K^i({\widehat M/x_1\widehat M})$ is Cohen-Macaulay of dimension $i$ for all $i\leq d-1.$  As $x_1$ is an $M$-sequential element, it follows by Lemma \ref{L:2bbs}(a) that  $x_1$ is a $K^i({\widehat M})$-regular element for all $i\geq 1$ with $K^i({\widehat M}))\neq 0$. 

 Since $x_1$ is an $M$-sequential element, $x_1$ is an f-element, see Remark \ref{R:1}(b). Therefore, 
 we have the exact sequences 
$$0\rightarrow H^i_{\frak m}(M)/x_1H^i_{\frak m}(M)\rightarrow H^i_{\frak m}(M/x_1M)\rightarrow (0:_{H^{i+1}_{\frak m}(M)}x_1)\rightarrow 0$$
for all $i\geq 0,$ see \cite[Lemma 2.4]{DN}.  So, we get the induced exact sequences
$$0\rightarrow K^{i+1}({\widehat M})/x_1K^{i+1}({\widehat M})\rightarrow K^i({\widehat M/x_1\widehat M})\rightarrow (0:_{K^i({\widehat M})}x_1)\rightarrow 0$$
for all $i\geq 0.$ Let $i\geq 1.$ As $(0:_{K^i({\widehat M})}x_1)=0$, we have $K^{i+1}({\widehat M})/x_1K^{i+1}({\widehat M})\cong K^i({\widehat M/x_1\widehat M})$. Therefore, if $K^i({\widehat M/x_1\widehat M})=0$ then $K^{i+1}({\widehat M})=0$ by Nakayama Lemma; if $K^i({\widehat M/x_1\widehat M})$ is Cohen-Macaulay of dimension $i$, then $K^{i+1}({\widehat M}) \neq 0,$ hence $x_1$ is $K^{i+1}({\widehat M})$-regular and hence $K^{i+1}({\widehat M})$ is Cohen-Macaulay of dimension $i+1$. Note that if $K^1({\widehat M})\neq 0$ then $x_1$ is $K^1({\widehat M})$-regular and hence $K^1({\widehat M})$ is Cohen-Macaulay of dimension $1.$ Therefore, $\widehat M$ is sequentially Cohen-Macaulay by Lemma \ref{L:1d}(a).  

The implication (a)$\Rightarrow$(c) follows by Lemma \ref{L:2bs} and by the above proof of (a)$\Rightarrow$(e).
 
 The rest statements (c)$\Rightarrow$(b) and (b)$\Rightarrow$(d) are obvious.
\end{proof}

\begin{remark} {\rm For any integer $d\geq 2$, there exists a Noetherian local ring of dimension $d$  such that it admits a sequential sequence s.o.p, but it is not a quotient of a Cohen-Macaulay ring (see Example \ref{E:1} for $d\geq 3$). In case where $d=2,$ let  $R$ be the  $2$-dimensional local domain  such that $\widehat R$ is sequentially Cohen-Macaulay but  $R$ is not sequentially Cohen-Macaulay (such a domain exists, see  \cite[Example 6.1]{Sch}). Let $a,b$ be a s.o.p of $R$. It is clear that  $a, b$ is an f-sequence s.o.p of $\widehat R$. So, $a, b$ is an $\widehat R$-sequential sequence by Theorem \ref{T:1}. Hence $a, b$ is an $R$-sequential sequence by Lemma \ref{L:2bbs}(b). Since $R$ is not sequentially Cohen-Macaulay, it follows by Theorem \ref{T:1} that $R$ is not a quotient of a Cohen-Macaulay ring.}
\end{remark}

\begin{lemma}  \label{L:2d} Suppose that $M$ be sequentially generalized Cohen-Macaulay. If $x \in \frak m$ is a generalized regular element of $M$ then $M/xM$ is sequentially generalized Cohen-Macaulay. 
\end{lemma}

\begin{proof}  If $d\leq 2$ then $\dim_R(M/xM)\leq 1$ and hence $M/xM$ is sequentially generalized Cohen-Macaulay. So, we assume that $d\geq 3.$ Denote by $M'$ the largest submodule of $M$ of dimension at most $1.$ Note  that $M/M'$ is sequentially generalized Cohen-Macaulay. Since  $\Ass_R(M/M')=\{\frak p \in \Ass_R(M)\mid \dim (R/\frak p) \geq 2\}$ by Lemma \ref{L:1c}(b) and $x$ is a generalized regular element of $M$, it follows that  $x$ is an $M/M'$-regular element. Therefore we get by \cite[Lemma 4.1(ii)]{NDC} that $(M/M')\big/x(M/M')$ is sequentially generalized Cohen-Macaulay. Set $\overline M:=M/xM$, $\overline M':=(M'+xM)/xM$.   Note that $\overline M/\overline M'\cong (M/M')\big/x(M/M')$. Therefore,  $\overline M/\overline M'$ is sequentially generalized Cohen-Macaulay.
Let  $$H^0_{\frak m}(\overline M/\overline M')=L_k/\overline M' \subset \ldots \subset L_1/\overline M'\subset L_0/\overline M'=\overline M/\overline M'$$  be the dimension filtration of $\overline M/\overline M',$ where $L_i$ is a submodule of $\overline M$ containing $\overline M'.$   Then $L_{i-1}/L_i$ is generalized Cohen-Macaulay for all integers $1\leq i\leq k.$ Set $H:=H^0_{\frak m}(M/xM).$ Consider the filtration of $\overline M=M/xM$
$$H\subseteq L_k+H \subset \ldots \subset L_1+ H \subset L_0=M/xM.$$ Because  $\dim_R(L_i/\overline M')<\dim_R(L_{i-1}/\overline M')$, we get  $\dim_R(L_i+H) <\dim_R(L_{i-1}+H)$ for any integer   $i \leq k$.  
On the other hand, from the  exact sequence
$$0\to L_{i-1}/L_i\to (L_{i-1}+H)/L_i \to (L_{i-1}+H)/L_{i-1}\to 0$$
with note that $(L_{i-1}+H)/L_{i-1}$ is of finite length and $L_{i-1}/L_i$ is generalized Cohen-Macaulay, it follows that $(L_{i-1}+H)/L_i $ is generalized Cohen-Macaulay for all $i \leq k$. So, from the exact sequence
$$0\to (L_{i}+H)/L_i\to (L_{i-1}+H)/L_i \to (L_{i-1}+H)/(L_{i}+H)\to 0$$
with notice that $(L_{i}+H)/L_i$ is of finite length, we have $(L_{i-1}+H)/(L_{i}+H)$ is generalized Cohen-Macaulay for all $i\leq k$. From the exact sequence $0\to \overline M'\to L_k\to L_k/\overline M'\to 0 $ with notice that  $L_k/\overline M'= H^0_{\frak m}(\overline M/\overline M')$ and $\dim_R(\overline M')\leq 1$, we have $\dim_R(L_k)\leq 1.$  Therefore, $(L_{k}+H)/H$ is generalized Cohen-Macaulay. Now, we get by \cite[Lemma 4.4(iii)]{CN} that $M/xM$ is sequentially generalized Cohen-Macaulay. 
\end{proof}

For  an $M$-regular $x$, if  $M$ is sequentially Cohen-Macaulay  then so is $M/xM$ (see  \cite[Theorem 4.7]{Sch}); however the converse statement is not true, see the counterexample in  \cite[Example 3]{CSSS}.  Similarly, the converse statement of Lemma \ref{L:2d} does not hold. For example, let $R=S/I$, where $S=k[[x,y,z,t,w]]$ be the formal power series ring of $5$ invariants over a field $k$ and $I=(x,y)\cap (z,t)$. Then $w$ is an $R$-regular  element, $R/wR$ is generalized Cohen-Macaulay, but $R$ is not sequentially generalized Cohen-Macaulay.  

Now, we are ready to prove Theorem \ref{T:2a}. 

\begin{theorem} \label{T:2}  The following statements are equivalent:
\begin{itemize}
\item[\rm{(a)}] $M$ is sequentially generalized Cohen-Macaulay.
\item[\rm{(b)}] $R/\Ann_R(M)$ is a quotient of a Cohen-Macaulay local ring and each generalized regular s.o.p of $M$ is an $M$-sequential f-sequence. 
\end{itemize}
\end{theorem}

\begin{proof} {\rm (a)} $\Rightarrow$ {\rm (b)}. We get by assumption (a) and by Lemma \ref{L:2bs} that $R/\Ann_R(M)$ is a quotient of a Cohen-Macaulay local ring. Let $(x_1, \ldots , x_d)$ be a generalized regular s.o.p of $M$. We need to show that  $(x_1, \ldots , x_d)$  is an $M$-sequential f-sequence. We prove this by induction on $d$. The case $d \leq 2 $ is obvious by Lemma \ref{L:1b}(c). Let $d\geq 3$. Suppose in contrary that  $x_1$ is not an $M$-sequential f-element. By Lemma \ref{L:2bbs}(a), $x_1\in\frak P$ for some $\frak P\in \Ass_{\widehat R}(K^i({\widehat M})) \setminus \{\frak m \widehat R\}$ with $i\geq 2.$  Hence $\dim_{\widehat R}(K^i({\widehat M}))>0.$ Note that $\widehat M$ is sequentially generalized Cohen-Macaulay by assumption (a) and by Lemma \ref{L:2b}(b).  Therefore we get by Lemma \ref{L:1d}(b) that  $K^i({\widehat M})$ is generalized Cohen-Macaulay of dimension $i$. Hence $\dim (\widehat R/\frak P)=i.$  It follows that   $\frak P\in\Ass_{\widehat R}(\widehat M)$. Hence $\frak p=\frak P\cap R\in\Ass_RM$ by \cite[Lemma 3.4]{Nh}. We note that $\dim (R/\frak p) \geq \dim (\widehat R/\frak P)=i \geq 2$.  Since $x_1$ is a generalized regular element of $M$ and $x_1 \in \frak p$, we have $\dim(R/\frak p) \leq 1.$ This gives  a contraction. So, $x_1$ is an $M$-sequential f-element.   Since $x_1$ is a generalized regular element of $M$, we get by assumption (a) and by Lemma \ref{L:2d} that $M/x_1M$ is sequentially generalized Cohen-Macaulay. Note that $x_2, \ldots, x_d$ is a generalized regular s.o.p of $M/x_1M$. So $x_2, \ldots, x_d$ is an $M/x_1M$-sequential f-sequence by induction. Thus, $x_1, \ldots, x_d$ is an $M$-sequential f-sequence.

{\rm (b)} $\Rightarrow$ {\rm (a)}. Since $R/\Ann_R(M)$ is a quotient of a Cohen-Macaulay local ring, $R/\frak p$ is unmixed for all $\frak p\in\Ass_R(M)$. So, Lemma \ref{L:2b}(b), it is enough to prove that $\widehat M$ is sequentially  generalized Cohen-Maccaulay. We prove this by induction on $d.$ The case  $d \leq 1$ is obvious. Let $d \geq 2$. Let $x_1 \in \frak m$ be a generalized regular element of $M$. Then $x_1$ is a part of a generalized regular sequence s.o.p of $M.$ Therefore, we get by assumption (b) that $x_1$ is an $M$-sequential f-element. Since $x_1$ is a generalized regular element of $M$,  we have the exact sequences
$$0\rightarrow H^j_{\frak m}(M)/x_1H^j_{\frak m}(M)\rightarrow H^j_{\frak m}(M/x_1M)\rightarrow (0:_{H^{j+1}_{\frak m}(M)}x_1)\rightarrow 0$$
for all integers $j\geq 1,$ see \cite[Lemma 2.4]{DN}.  So,  we have the induced exact sequences
$$\ \ \ \ \ \ \ \ \ \ \ \ \ \ \ \ \ \ \ \ 0\rightarrow K^{j+1}({\widehat M})/x_1K^{j+1}({\widehat M})\rightarrow K^j({\widehat M/x_1\widehat M})\rightarrow (0:_{K^j({\widehat M})}x_1)\rightarrow 0 \ \ \ \ \ \ \ \ \ \ \ \ \ \ \ \ \ \ \ \ (1)$$
for all $j\geq 1.$ Note that $M/x_1M$ satisfies the induction hypothesis. So, $\widehat M/x_1\widehat M$ is sequentially generalized Cohen-Macaulay. By Lemma \ref{L:1d}(b), either  $K^j({\widehat M/x_1 \widehat M})$ is  generalized Cohen-Macaulay of dimension $j$ or $\ell_{\widehat R}(K^j({\widehat M/x_1 \widehat M}))< \infty$  for all $j\leq d-1.$ First, we use Lemma \ref{L:1g} to prove the following claim. 

 \noindent {\it Claim:  The sequential polynomial type ${\rm sp}(M)$ of $M$ is at most $1.$}

 In fact, since $R/\Ann_R(M)$ is a quotient of a Cohen-Macaulay local ring by assumption (b), $R/\frak p$ is unmixed for all $\frak p\in\Ass_R(M)$. Hence ${\rm sp}(M)={\rm sp}(\widehat M)$ by \cite[Theorem 3.5]{NDC}. Therefore, we need only to prove that ${\rm sp}(\widehat M)\leq 1.$ Set $D(\widehat M)=\{\dim(\widehat R/\frak Q)\mid \frak Q\in\Ass_{\widehat R}(\widehat M)\}.$ By Lemma \ref{L:1g}, it is enough to prove that $\dim_{\widehat R}\big(K^i({\widehat M})\big)\leq 1$ for all $i\notin D(\widehat M)$ and  ${\rm p}(K^i({\widehat M}))\leq 1$ for all $i\in D(\widehat M)$. If $i\leq 1$ then  ${\rm p}(K^i({\widehat M}))\leq \dim_{\widehat R}\big(K^i({\widehat M})\big)\leq 1$ by Lemma \ref{L:1e}(a). Let $i=2.$ It follows by  \cite[10.2.20, 11.2.6]{BS} and Lemma \ref{L:1b}(a) that $\dim_{\widehat R}\big(K^2({\widehat M})\big)\leq 1$ if $2\notin D(\widehat M)$ and $\dim_{\widehat R}\big(K^2({\widehat M})\big)=2$ if $2\in D(\widehat M).$ Therefore, ${\rm p}(K^2({\widehat M}))\leq 1$ by Lemma \ref{L:1e}(a).

 Let  $i\geq 3$. If $\ell_{\widehat R}\big(K^{i-1}({\widehat M/x_1 \widehat M})\big)< \infty$ then $\ell_{\widehat R}\big(K^i({\widehat M})/x_1K^i({\widehat M})\big)<\infty$ by exact sequences (1), hence ${\rm p}(K^i({\widehat M}))\leq \dim_{\widehat R}\big(K^i({\widehat M})\big)\leq 1.$  So, we assume that $K^{i-1}({\widehat M/x_1 \widehat M})$ is generalized Cohen-Macaulay of dimension $i-1$. Since $x_1$ is an $M$-sequential f-element,  we get by Lemma \ref{L:2bbs}(a) that $x_1$ is an f-element of $K^j({\widehat M})$ for all $j\geq 2$. Hence $\ell_{\widehat R}(0:_{K^j({\widehat M})}x_1)<\infty$ for all $j\geq 2.$    As $i-1\geq 2,$ it follows by the exact sequence (1) at $j=i-1$ that $K^i({\widehat M})/x_1K^i({\widehat M})$ is generalized Cohen-Macaulay of dimension $i-1$. Therefore we have $\dim_{\widehat R}(K^i({\widehat M}))\geq i-1>0.$  Since $x_1$ is an f-element of $K^i({\widehat M})$, it follows that $x_1$ is a parameter element of $K^i({\widehat M})$ and hence $\dim_{\widehat R}(K^i({\widehat M}))=i.$ So, there exists $\frak P\in\Ass_{\widehat R}(K^i({\widehat M}))$ such that $\dim(\widehat R/\frak P)=i.$ Therefore,  $\frak P\in\Ass_{\widehat R}(\widehat M)$. Hence $i\in D(\widehat M).$ Set $K^i:=K^i(\widehat M).$ Since $x_1$ is an f-element of $K^i$, we have the exact sequences 
$$0\rightarrow H^j_{\frak m\widehat R}(K^i)/x_1H^j_{\frak m\widehat R}(K^i) \rightarrow H^j_{\frak m\widehat R}(K^i/x_1K^i) \rightarrow (0:_{H^{j+1}_{\frak m\widehat R}(K^i)}x_1) \rightarrow 0$$
for all $0\leq j\leq i-1,$ see \cite[Lemma 2.4]{DN}. Since $K^i/x_1K^i$ is generalized Cohen-Macaulay of dimension $i-1$, we have $\ell_{\widehat R}\big(0:_{H^j_{\frak m\widehat R}(K^i)}x_1\big)<\infty$ for all $1\leq j\leq i-1.$ It follows by Local Duality Theorem \cite[11.2.6]{BS} that $\dim_{\widehat R}\big(\widehat R/\Ann_{\widehat R}H^j_{\frak m\widehat R}(K^i)\big)\leq 1$ for all $1\leq j\leq i-1.$ Therefore, ${\rm p}(K^i(\widehat M))\leq 1$ by Lemma \ref{L:1e}(a). Thus, ${\rm sp}(M)={\rm sp}(\widehat M)\leq 1$ by Lemma \ref{L:1g}, the claim follows.

Now we continue the proof of (b)$\Rightarrow$(a). Suppose in contrary that $\widehat M$ is not sequentially gneralized Cohen-Macaulay. Then ${\rm sp}(M)={\rm sp}(\widehat M)=1$ by Claim. Moreover, there exists  by Lemma \ref{L:1d}(b) an integer $i$ such that  $\dim_{\widehat R}(K^i({\widehat M}))>0$ and $K^i({\widehat M})$ is not generalized Cohen-Macaulay of dimension $i$. Note that $K^1({\widehat M})$ is either of finite length or is generalized Cohen-Macaulay of dimension $1$. Therefore, $i\geq 2.$ We consider two cases.

$\bullet${\it Case 1:}  $\ell_{\widehat R}(K^{i-1}(\widehat M/x_1 \widehat M))< \infty$

We get by the  exact sequences (1) at $j=i-1$ that $\ell_{\widehat R}(K^i(\widehat M)/x_1K^i(\widehat M))<\infty$. Since $\dim_{\widehat R}(K^i({\widehat M}))>0$, we have  $\dim_{\widehat R}\big(K^i({\widehat M}))=1.$ Let $\frak P\in\Ass_{\widehat R}\big(K^i({\widehat M}))$ such that $\dim (\widehat R/\frak P)=1.$ We have by  \cite[10.2.20, 11.2.6]{BS} that $\frak P\in\Att_{\widehat R}H^i_{\frak m}(M)$. Set $\frak p=\frak P\cap R.$ Then $\frak p\in \Att_R H^i_{\frak m}(M)$ by Lemma \ref{L:1a}(b). As $R/\Ann_R(M)$ is a quotient of a Cohen-Macaulay local ring, $R/\frak p$ is unmixed. Moreover, we get by the shifted principle for attached primes of local cohomology modules under completion (see Lemma \ref{L:NQ}(a)) that $\frak P\in\Ass (\widehat R/\frak p\widehat R)$.  Therefore, $\dim (R/\frak p)=\dim (\widehat R/\frak P)=1.$ Since $d \geq 2$, there exists a s.o.p $y_1, \ldots , y_d$ of $M$ which is an $M$-generalized regular sequence such that $y_1\in\p$. Hence,  $y_1, \ldots , y_d$ is an $M$-sequential f-sequence by assumption (b). Since $i \geq 2$ and $\frak p\in \Att_R H^i_{\frak m}(M),$ it follows that $y_1\notin \frak p$. This gives a contradiction.

$\bullet${\it Case 2:}  $K^{i-1}(\widehat M/x_1 \widehat M)$ is generalized Cohen-Macaulay of dimension $i-1$

As in the proof of the claim, we have $i\in D(\widehat M)$, $\dim_{\widehat R}(K^i({\widehat M}))=i$ and ${\rm p}(K^i({\widehat M}))\leq 1.$ Therefore, $K^i({\widehat M})$ is not generalized Cohen-Macaulay. Hence ${\rm p}(K^i({\widehat M}))=1.$ So, it follows by Lemma \ref{L:1e}(b) that 
$$\max\{\dim_{\widehat R} \nCM (K^i({\widehat M})), \dim_{\widehat R}\big(U_{K^i({\widehat M})}(0)\big)\}=1,$$
where $\nCM (K^i({\widehat M}))$ is the non Cohen-Macaulay locus of $K^i({\widehat M})$ and $U_{K^i({\widehat M})}(0)$ is the largest submodule of $K^i({\widehat M})$ of dimension less than $i.$ 

Suppose that $\dim_{\widehat R}\big(U_{K^i({\widehat M})}(0)\big)=1.$ Then there exists $\frak Q\in \Ass_{\widehat R}(K^i({\widehat M}))$ such that $\dim(\widehat R/\frak Q)=1.$ Hence $\frak Q\in \Att_{\widehat R}H^i_{\frak m}(M)$ by \cite[10.2.20, 11.2.6]{BS}. Set $\frak q=\frak Q\cap R.$ Then $\frak q\in\Att_RH^i_{\frak m}(M)$ by Lemma \ref{L:1a}(b). As $R/\Ann_R(M)$ is a quotient of a Cohen-Macaulay local ring, $R/\frak q$ is unmixed. Moreover, we get by shifted principle for attached primes of local cohomology modules under completion (see Lemma \ref{L:NQ}(a)) that $\frak Q\in\Ass (\widehat R/\frak q\widehat R)$. Therefore, $\dim (R/\frak q)=\dim(\widehat R/\frak Q)=1.$ As $d \geq 2$, there exists a s.o.p $y_1, \ldots , y_d$ of $M$ which is a generalized regular sequence such that $y_1\in\q$. By assumption (b), $y_1, \ldots , y_d$ is a sequential f-sequence of $M$. As $i\geq 2$ and $\frak q\in\Att_RH^i_{\frak m}(M)$, we get $y_1\notin \frak q$. This gives a contradiction. 

Suppose that  $\dim_{\widehat R} \nCM (K^i({\widehat M}))=1.$ Let $\frak Q\in \nCM (K^i({\widehat M}))$ such that $\dim (\widehat R/\frak Q)=1.$ Then $K^{i-1}({\widehat M_{\frak Q})}$ is not Cohen-Macaulay. Hence $\widehat M_{\frak Q}$ is not sequentially Cohen-Macaulay by Lemma \ref{L:1d}(a). Let $\frak q=\frak Q\cap R.$ Then $\dim (R/\frak q)\geq 1.$ If $\dim (R/\frak q)\geq 2$ then we get by the claim and by \cite[Theorem 3.4(ii)]{NDC} that ${\rm sp}(M_{\frak q})\leq {\rm sp}(M)-\dim(R/\frak q)<0.$ It means that  $M_{\frak q}$ is sequentially Cohen-Macaulay. Since the natural map $R_{\frak q}\to \widehat R_{\frak Q}$ is a flat local homomorphism and the fiber $\widehat R_{\frak Q}/\frak q\widehat R_{\frak Q}$ is Cohen-Macaulay by assumption (b), it follows by \cite[Theorem 5]{TY} that $\widehat M_{\frak Q}$ is sequentially Cohen-Macaulay. This is impossibe. Therefore, we can assume that $\dim(R/\frak q)=1.$ Since $d \geq 2,$ there exists a s.o.p $y_1, \ldots , y_d$ of $M$ which is an $M$-generalized regular sequence such that $y_1, \ldots , y_{d-1}\in\frak q$. Hence,  $y_1, \ldots , y_d$ is an $M$-sequential f-sequence by assumption (b). Set $\dim_{R_{\frak q}}(M_{\frak q})=r.$ Then $r\leq d-1$. Since $y_1, \ldots , y_r$ is a generalized regular sequence of $M,$ it follows that the elements $y^*_1, \ldots , y^*_r\in \frak q R_{\frak q}$ form  an f-sequence s.o.p. of $M_{\frak q}$, where $y^*_j$ is the image of $y_j$ in $R_{\frak q}$. Moreover, since $R/\Ann_R(M)$ is a quotient of a Cohen-Macaulay local ring and $y_1, \ldots , y_r$ is an $M$-sequential f-sequence, it follows by the shifted principle for attached primes of local cohomology modules under localization (see Lemma \ref{L:NQ}(b)) that $y^*_1, \ldots , y^*_r$ is an $M_{\frak q}$-sequential sequence.  So, we have by Theorem \ref{T:1}(b)$\Rightarrow$(a) that $M_{\frak q}$ is sequentially Cohen-Macaulay. Hence  $\widehat M_{\frak Q}$ is sequentially Cohen-Macaulay by \cite[Theorem 5]{TY}. This gives a contradiction.
\end{proof}

As an application, we have the following new characterizations of Cohen-Macaulay modules and generalized Cohen-Macaulay modules in terms of sequential sequence and sequential f-sequence.

\begin{corollary} \label{C:3}  Suppose that $d\geq 2.$ Denote by $M'$ the largest submodule of $M$ of dimension at most $1$. Then 
\begin{itemize}
\item[\rm{(a)}]  $M$ is Cohen-Macaulay if and only if $R/ \Ann_R(M)$ is a quotient of a Cohen-Macaulay local ring, $H^0_{\frak m}(M)=0$ and each s.o.p of $M$ is an $M$-sequential sequence.
\item[\rm{(b)}]  $M$ is generalized Cohen-Macaulay if and only if  $R/ \Ann_R(M)$ is a quotient of a Cohen-Macaulay local ring, $M'=H^0_{\frak m}(M)$ and each s.o.p of $M$ is an $M$-sequential f-sequence.
\end{itemize}
\end{corollary}

\begin{proof} (a). If   $M$ is Cohen-Macaulay then it is clear that   $R/ \Ann_R(M)$ is a quotient of a Cohen-Macaulay local ring, $H^0_{\frak m}(M)=0$ and each s.o.p of $M$ is an $M$-sequential sequence. Conversely, we have by Theorem \ref{T:1} that $M$ is sequentially Cohen-Macaulay. Moreover, we get by Remark \ref{R:1}(b) that each s.o.p of $M$ is an f-sequence. Since $R/ \Ann_R(M)$ is a quotient of a Cohen-Macaulay local ring, $M$ is generalized Cohen-Macaulay. So,  $H^0_{\frak m}(M)\subset M$ is the dimension filtration of $M$. Hence $M$ is Cohen-Macaulay as $H^0_{\frak m}(M)=0$.

(b). If  $M$ is generalized Cohen-Macaulay then it is clear that   $R/ \Ann_R(M)$ is a quotient of a Cohen-Macaulay local ring, $M'=H^0_{\frak m}(M)$ and each s.o.p of $M$ is an $M$-sequential f-sequence. Conversely, we have by Theorem \ref{T:2} that $M$ is sequentially generalized Cohen-Macaulay. Moreover, we get by Remark \ref{R:1}(b) that each s.o.p of $M$ is a generalized regular sequence. For any $\frak p\in\Spec(R)$, if $1<\dim(R/\frak p)<d$ then we can choose a s.o.p. $x_1, \ldots , x_d$ of $M$ such that $x_1\in\frak p$. Then  $x_1, \ldots , x_d$ is a generalized regular sequence of $M$, therefore $\frak p\notin \Ass_R(M).$ So, $M'$ is the largest submodule of $M$ of dimension less than $d$. Since $M'=H^0_{\frak m}(M),$ the dimension filtration of $M$ is  $H^0_{\frak m}(M)\subset M$.  Hence $M/H^0_{\frak m}(M)$ is generalized Cohen-Macaulay and hence $M$ is generalized Cohen-Macaulay. 
\end{proof}

Finally, we give an example to clarify the results in Theorems \ref{T:1}, \ref{T:2}. The following example shows that for a given integer $d\geq 3$,  there exists a  Noetherian local domain $(R,\frak m)$ of dimension $d$ such that

a)  $R$ admids an $R$-sequential sequence s.o.p;

b) $R$ is not a quotient of a Cohen-Macaulay ring;

c) $\widehat {R}$ is sequentially Cohen-Macaulay; $R$ is not sequentially generalized Cohen-Macaulay.

\begin{example} \label{E:1} {\rm  Let $(V, \frak v)$ be the Noetherian local domain of dimension $3$ constructed by M. Brodmann and C. Rotthaus \cite[(15)]{BR} such that the completion $\widehat {V}$ of $V$ can be identified with $S/I$, where $S =\Bbb Q[[x, y, z, t]]$ is the formal power series ring of invariants $x, y, z, t$ over $\Bbb Q$ and $I = (x)\cap (y)\cap (x^2, y^2).$ Since $(x,y)$ is an embedded prime ideal of $\widehat {V}$, it follows that $R$ is not a quotient of a Cohen-Macaulay ring. By \cite[Example 3.7]{NDC}, we have $\widehat {V}$ is sequentially Cohen-Macaulay, $\depth (V)=2$ and ${\rm sp}(V)={\rm p}(V)=2$. In particular, $V$ is not sequentially generalized Cohen-Macaulay.  Hence $\Att_{\widehat V}(H^i_{\frak v}(V))=\emptyset$ for $i=0,1$ and $K^2(\widehat {V})$ is Cohen-Macaulay of dimension $2.$ So, we get by Lemma \ref{L:1b}(a) that $\Att_{\widehat {V}}(H^2_{\frak v}(V))=\{(x,y)\}.$ Hence, $\Att_{V}(H^2_{\frak v}(V))=\{0\}$ by Lemma \ref{L:1a}(b). Let $a_1, a_2, a_3$ be an f-sequence s.o.p. of $V.$ Then it is clear that  $a_1, a_2, a_3$ is an f-sequence s.o.p. of $\widehat {V}$. Since $\widehat {V}$ is sequentially Cohen-Macaulay, $a_1, a_2, a_3$ is an $\widehat {V}$-sequential sequence by Theorem \ref{T:1}. So, by Lemma \ref{L:2bbs}(b), $a_1, a_2, a_3$ is an $V$-sequential sequence.  

If $d=3$ then we take $R=V.$

Now, let $d>3.$ Set $d=3+r$. Let $R:= V[[x_1, \ldots , x_r]]$, where
$x_1,\ldots , x_r$ are independent indeterminates over $V$. Let $\frak m$ be the unique maximal ideal of $R$.  Then $R$ is not a quotient of a Cohen-Macaulay ring since $\widehat R$ has an embedded prime ideal.   By \cite[Example 3.7]{NDC},  $\widehat {R}$ is a sequentially Cohen-Macaulay local domain of dimension $d$ with  $\depth (R)=d-1$ and ${\rm sp}(R) = {\rm p}(R)=d-1$. By the same arguments as in the above, $\Att_R(H^{d-1}_{\frak m}(R))=\{0\}$ and each f-sequence s.o.p. of $R$ is an $R$-sequential sequence.}
\end{example}

\noindent {\bf Acknowledgment}

We wish to express our gratitude to the referee for many useful suggestions.

\end{document}